\documentclass[12pt]{article}
\usepackage{float}
\usepackage{longtable}
\usepackage{amsmath}
\usepackage{amssymb}
\usepackage{theorem}
\usepackage{graphicx} 
\newtheorem{theorem}{Theorem}
\newtheorem{lemma}{Lemma}
\newtheorem{proposition}[theorem]{Proposition}
\newtheorem{example}{Example}

\newcommand{\cay}{\mathrm{Cay}}

\author{\v Stefan Gy\"urki\footnote{e-mail address: {\tt gyurki@savbb.sk}} \\[0.5cm]
{\small Matej Bel University, 974 11 Bansk\'a Bystrica, Slovak Republic}} 
\title{New infinite families \\ of directed strongly regular graphs \\
via equitable partitions}
\date{\today}

\begin{document}

\maketitle

\begin{abstract}
In this paper we introduce a construction of directed strongly regular graphs
from smaller ones using equitable partitions. Each equitable partition of a single
DSRG satisfying several conditions leads to an infinite family of directed strongly 
regular graphs. We construct in this way dozens of infinite families. 
For order at most 110, we confirm the existence of DSRGs for 30 previously open 
parameter sets.  
\end{abstract}

\section{Introduction}
The main subject of this paper is directed strongly regular graph
(briefly DSRG), a possible generalisation of the well-known (undirected) 
strongly regular graphs for the directed case introduced by Duval in \cite{du}. 
The undirected version plays a central role in Algebraic Graph Theory, 
while the directed version received less attention. 

We show how to construct larger DSRGs from smaller ones 
with the aid of suitably defined equitable partitions for directed graphs. 
According to our best knowledge, equitable partitions for directed graphs
have not been defined in exact terms, though similar concept can be found
for association schemes in \cite{go} with an application for (undirected) 
strongly regular graphs in \cite{hi}.

The main results of the introduced idea involve constructions of many infinite families
of DSRGs. Several families have new parameter sets, others 
generalise and connect constructions which have been previously known
from various ideas.

For the small parameter sets there is an evidence of the current state-of-the-art 
on the webpage of A.Brouwer and S.Hobart (see \cite{ab}). According to their
catalogue of parameter sets with order at most 110, 
we have confirmed the existence of DSRGs for 30 open parameter sets. 

The paper is organised as follows: In Section 2 the necessary and most significant 
terms are repeated.
In Section 3 a brief summary is given on the used computer algebra packages.  
In Section 4 we define a new construction called $\pi$-join. It serves as a key concept
for our project. We derive the necessary and sufficient conditions, when this construction
gives us a DSRG if the basic graph is a smaller DSRG. It turns out that the crucial 
moment is finding a certain partition in the basic graph. In Section 5 we present 
good partitions for several well-known families of DSRGs and in Section 6 for DSRGs with 
some sporadic parameter sets. 
 
\section{Preliminaries}

In this section we present the basic and most significant concepts used in this paper. 

\subsection{General concepts in algebraic graph theory}

All graphs considered in this paper are finite, simple, directed or undirected graphs, 
thus no loops or multiple edges are allowed. If there is a dart from vertex $x$ 
to vertex $y$ in a directed graph (digraph), then we say that $y$ is an \emph{outneighbour}
of $x$, and $x$ is an \emph{inneighbour} of $y$. The number of outneighbours 
(inneighbours) of a vertex $x$ is called \emph{outdegree (indegree)} of $x$. 
If the in- and outdegrees of all vertices in a digraph are equal to $k$, then we say that 
the digraph is \emph{$k$-regular}, or \emph{regular of degree (valency) $k$}.

Let $H$ be a group. Suppose that $X\subseteq H$ and $e\notin X$, 
where $e$ is the identity element of the group $H$. 
Then the digraph $\Gamma=\cay(H,X)$ with vertex set $H$ and
dart set $\{(x,y): x,y\in H, yx^{-1}\in X\}$ is called the {\it
Cayley digraph over $H$ with respect to~$X$.}

Let $\Gamma=(V,E)$ be a simple undirected graph. Assume that $\pi=\{C_1,\ldots,C_r\}$
is a partition of $V$. It is called \emph{equitable} if for all $1\leq i,j\leq r$ the
amount of neighbours of a vertex $u\in C_i$ in the set $C_j$ is equal to a constant
$q_{ij}$ which does not depend on the choice of $u\in C_i$. Thus, to an equitable
partition we can associate a matrix $Q=(q_{i,j})_{1\leq i,j\leq r}$ which is called
\emph{quotient matrix}. An equitable partition $\pi$ is called \emph{homogeneous of degree $d$},
if each cell has size $d$. Equitable partitions of graphs can be defined also using
their adjacency matrix. We say that a partition $\pi$ of a matrix $A$ 
is \emph{column equitable}, if there exists such an $r\times r$ matrix $Q$, 
usually called \emph{quotient matrix}, for which $A^TH=HQ$, where
$H$ is an $n\times r$ matrix, the \emph{characteristic matrix} of $\pi$,
whose columns are the characteristic vectors of $\pi$. In other words, 
$h_{ij}=\begin{cases} 1 & \textrm{ if } v_i\in C_j;\\
0 & \textrm{ otherwise.} \end{cases}$ 

Similarly, a partition $\pi$ of a matrix $A$ is \emph{row equitable},
if $\pi$ is column equitable for $A^T$. Clearly, if $\pi$ is both row and column
equitable for an adjacency matrix $A$, then $\pi$ defines an equitable partition
of the corresponding graph.  

Following this idea one can define a kind of equitable partitions also for directed graphs.
We will say that $\pi$ is a \emph{column equitable partition} of a directed graph $\Gamma$, 
if it is a column equitable partition of its adjacency matrix $A(\Gamma)$. 

Let $A$ be an adjacency matrix of a (di)graph $\Gamma$ of order $n$, 
and $\pi$ a partition of $V(\Gamma)$ into $a$ cells $C_1,C_2,\ldots,C_a$.
We say that adjacency matrix $A$ \emph{respects} the partition $\pi$, if the 
corresponding ordering $(v_1,v_2,\ldots,v_n)$ of the vertices satisfies: 
$i\leq j \iff x\leq y$ for all vertices $v_i\in C_x$, $v_j\in C_y$.

We assume that the reader is familiar with the terms like \emph{adjacency matrix}, 
\emph{Kronecker product of matrices}, \emph{group of automorphisms of a graph}, 
\emph{1-factor (perfect matching)}, \emph{clique, coclique} and 
\emph{line graphs}. Otherwise, we refer the reader to \cite{cl}. 

We will use also terms like \emph{blocks of imprimitivity, rank of a permutation group}
or \emph{transitivity} from permutation group theory. For definitions see \cite{dm}. 

\subsection{Strongly regular graphs}

An undirected graph $\Gamma$ is called \emph{regular} of valency $k$, or \emph{$k$-regular},
if each vertex is incident to the same number $k$ of edges.
In the terms of adjacency matrix of a graph we can equivalently say that 
$\Gamma$ is $k$-regular if and only if for its adjacency matrix $A=A(\Gamma)$ the equation
$AJ=JA=kJ$ holds, where $J$ is the all-one matrix.  
A simple regular graph with valency $k$ is said to be \emph{strongly regular}
(SRG, for short) if there exist integers $\lambda$ and $\mu$ such that for each edge $\{u,v\}$ the number
of common neighbours of $u$ and $v$ is exactly $\lambda$; while for each non-edge
$\{u,v\}$ the number of common neighbours of $u$ and $v$ is equal to $\mu$. 
Previous condition can be equivalently rewritten into the equation $A^2=kI+\lambda A+\mu(J-I-A)$
using the adjacency matrix of $\Gamma$.
The quadruple $(n,k,\lambda,\mu)$ is called the \emph{parameter set} of an SRG $\Gamma$. 

Strongly regular graphs play a central importance in Algebraic Graph Theory. There are
many dozens of papers with interesting constructions and results. The most significant
ones altogether with crucial information on their properties are collected on 
the webpage of A. Brouwer~\cite{ab}. A database of small SRGs is located on 
the webpage of E. Spence~\cite{ts}.  

\subsection{Directed strongly regular graphs}

A possible generalisation of the notion of SRGs for directed graphs was given by Duval~\cite{du}.
While the family of SRGs has been well-studied, the directed version has received less attention.
 
A \emph{directed strongly regular graph} (DSRG) with parameters $(n,k,t,\lambda,\mu)$ is a
regular directed graph on $n$ vertices with valency $k$, such that every vertex is incident with~$t$
undirected edges, and the number of directed paths of length 2 directed from a vertex $x$ 
to another vertex $y$ is
$\lambda$, if there is an arc from $x$ to $y$, and $\mu$ otherwise. 
In particular, a DSRG with $t=k$ is an SRG, and a DSRG with $t=0$ is a doubly regular tournament.

\medskip

\noindent{\bf Remark.}
In this paper we prefer to use the ordering $(n,k,t,\lambda,\mu)$ for the 5-tuple of parameters
of a DSRG, however in several other papers the ordering $(n,k,\mu,\lambda,t)$ is used. 

\medskip

The adjacency matrix $A=A(\Gamma)$ of a DSRG with parameters $(n,k,t,\lambda,\mu)$, satisfies 

\begin{equation}
 AJ=JA=kJ\quad\textrm{ and }\quad A^2=tI+\lambda A+\mu(J-I-A). \label{am}
\end{equation}
 
Oppositely, if an $n\times n$ zero-one matrix $A$ with zeros on the diagonal satisfies (\ref{am})
for some $k,t,\lambda$ and $\mu$, then it is an adjacency matrix of 
a DSRG with parameters $(n,k,t,\lambda,\mu)$.

\begin{proposition}[\cite{du}]
If $\Gamma$ is a DSRG with parameter set $(n,k,t,\lambda,\mu)$ and adjacency matrix $A$, 
then the complementary graph $\bar\Gamma$ is a DSRG with parameter set $(n,\bar k,\bar t,
\bar\lambda,\bar\mu)$ with adjacency matrix $\bar A=J-I-A$, where 
\begin{eqnarray*}
\bar k &=& n-k-1 \\
\bar t &=& n-2k+t-1 \\
\bar \lambda & =& n-2k+\mu-2 \\
\bar \mu &=& n-2k+\lambda.
\end{eqnarray*}
\label{prop2}
\end{proposition}

For a directed graph $\Gamma$ let $\Gamma^T$ denote the digraph obtained by reversing
all the darts in $\Gamma$. Then $\Gamma^T$ is called the \emph{reverse} of $\Gamma$. In other words,
if $A$ is the adjacency matrix of~$\Gamma$, then $A^T$ is the adjacency matrix of $\Gamma^T$.

The following proposition was observed by Pech, and presented in \cite{km}:

\begin{proposition}[\cite{km}]
Let $\Gamma$ be a DSRG. Then the graph $\Gamma^T$ is a DSRG with the same parameter set. 
\label{prop1}
\end{proposition}

We say that two DSRGs $\Gamma_1$ and $\Gamma_2$ are \emph{equivalent}, if $\Gamma_1\cong\Gamma_2$, 
or $\Gamma_1\cong\Gamma_2^T$, or $\Gamma_1\cong\bar\Gamma_2$, or $\Gamma_1\cong\bar\Gamma_2^T$;
otherwise they are called \emph{non-equivalent}. (In other words, $\Gamma_1$ is equivalent to 
$\Gamma_2$ if and only if $\Gamma_1$ is isomorphic to $\Gamma_2$ or to a graph, obtained from
$\Gamma_2$ via reverse and complementation.) 
From our point of view the interesting DSRGs are those which are non-equivalent. 

Under \emph{eigenvalues} of a DSRG we mean the eigenvalues of its adjacency matrix. 
It is known, that there are precisely three different eigenvalues, all of them are
integers, and their values together with their multiplicities are uniquely 
determined by the parameter set. 
More explicitly, a DSRG with parameter set $(n,k,t,\lambda,\mu)$ has eigenvalues:

\[
\theta_0=k,\quad \theta_{1,2}=\frac12\left(\lambda-\mu\pm\sqrt{(\mu-\lambda)^2+
4(t-\mu)}\right),
\] 
with respective multiplicities:
\[
m_0=1,\quad m_1=\frac{k+\theta_1(n-1)}{\theta_2-\theta_1},\quad
m_2=\frac{k+\theta_2(n-1)}{\theta_1-\theta_2}.
\]

The parameters $n,k,t,\lambda,\mu$ are not independent. Relations, which have to be satisfied 
for such parameter sets are usually called \emph{feasibility conditions}. 
Most important, and, in a sense, basic conditions are the following (for their proof see \cite{du}):

\begin{equation}
k(k+\mu-\lambda)=t+(n-1)\mu. \label{c1}
\end{equation}

There exists a positive integer $d$ such that:
\begin{eqnarray}
d^2&=&(\mu-\lambda)^2+4(t-\mu) \label{c2}\\
d &\mid& (2k-(\mu-\lambda)(n-1)) \label{c3}\\
n-1 &\equiv& \frac{2k-(\mu-\lambda)(n-1)}{d}  \pmod 2 \label{c4}\\
n-1 &\geq &\left|\frac{2k-(\mu-\lambda)(n-1)}{d}\right|.\label{c5}
\end{eqnarray}

Further:
\[
\begin{array}{rcccl}
0 & \leq & \lambda < t & < & k \\
0 & < & \mu \leq t & < & k \\
-2(k-t-1) & \leq & \mu  -  \lambda & \leq & 2(k-t).  
\end{array}
\]

We have to mention that for a feasible parameter set 
it is not guaranteed that a DSRG does exist. 
A feasible parameter set for which at least one DSRG $\Gamma$ exists is called \emph{realizable},
otherwise \emph{non-realizable}.
The smallest example of a feasible, but non-realizable parameter set is $(14,5,4,1,2)$. 
The non-existence of such graph was shown in \cite{km}.

There are known several general constructions of DSRGs. The most common 
constructions are based on algebraic or combinatorial approaches, or even 
on their combination. They include approaches using block matrices \cite{ad,du}, 
Kronecker product of matrices \cite{du}, finite geometries \cite{fm,fp,gh,kp}, 
combinatorial block designs \cite{fp,gy}, finite incidence structures \cite{bo,os}, 
coherent algebras \cite{fp,km}, association schemes \cite{gy}, 
regular tournaments \cite{jg}, partial sum families \cite{ar,ma}, 
Cayley graphs \cite{di,gy,hs,jg,km} or computer algebra experimentation
\cite{gy,jn}. In this paper we extend this list by a construction which
is on the border of matrix construction, equitable partitions 
and computer algebra experimentation.  

\section{Computer algebra tools}

This project has been supported by an essential amount of results which 
have been obtained with the aid of a computer. All such computations have
been run in {\sf GAP} (Groups, Algorithms, Programming \cite{ga}) jointly with 
using its share packages: 
{\sf GRAPE} -- a package for computation with graphs \cite{lh} together with 
{\sf nauty} for testing automorphisms of graphs \cite{na}; 
{\sf COCO II} -- an unreleased package for computation with coherent configurations
and association schemes, written by Reichard \cite{sr}; and the package {\sf SetOrbit}
for finding representatives of orbits of group actions on sets of various size, 
written by Pech and Reichard \cite{so,pr}. 

\section{$\pi$-join construction}
\label{spi}

In this section we introduce a new construction, which serves as a basic idea for
constructing new graphs from old ones. We will investigate when this construction
gives us a larger DSRG from smaller ones and derive its most important properties.
Starting from this section it often happens that we do not need to distinguish
between directed and undirected strongly regular graphs. In such a case we
use abbreviation (D)SRG. 

\subsection{Description of the construction}

Let us consider a directed or an undirected strongly regular graph $\Gamma$
of order $n$. Let us call it \emph{basic graph}. Suppose that $\pi=\{C_1,C_2,\ldots,C_a\}$ 
is a homogeneous partition of its vertex set into $a$ cells of size $b$, where $ab=n$. 
Consider a positive integer $j$ and define a digraph $\Gamma^j_{\pi}$ as following.
Create $ja+1$ disjoint copies $\Gamma_0,\Gamma_1,\ldots,\Gamma_{ja}$ of the graph 
$\Gamma$, i.e. for its vertex set we have
$$V(\Gamma^j_{\pi})=V(\Gamma)\times \{0,1,\ldots,ja\}.$$
In each of the copies we can distinguish the images of the cells of the partition
$\pi$. Let us denote $C_s^r$ the image of the cell $C_s$ in the $r$-th copy of 
$\Gamma$. We define the adjacency in $\Gamma^j_{\pi}$ between vertices 
by giving the set of outneighbours of vertices. 
For fixed $r$ there is a dart from vertex $(u,r)$ to vertex $(v,r)$ 
whenever there is a dart from $u$ to $v$ in $\Gamma$.
Further, let $(u,r)\in C_s^r$ for some $s$ and $r$. There is a dart from $(u,r)$ 
to all of the vertices in cells $C_1^{r+1}$, $C_1^{r+2}, \ldots,
C_1^{r+j}, C_2^{r+j+1}, C_2^{r+j+2}, \ldots, C_2^{r+2j}, \ldots, 
C_a^{r+(a-1)j+1}, \ldots, C_a^{r+aj}$, where the superscripts are
taken modulo $aj+1$. Let us call this construction \emph{$\pi$-join of $\Gamma$ in 
power $j$,} or just shortly \emph{$\pi$-join} when $\Gamma$ and/or $j$ are either clear
from the context or irrelevant.

According to our best knowledge this construction is new. However, 
it can be regarded as a mixture of the lexicographic
product of graphs \cite{sa} and the X-join of graphs \cite{sb}. 
In particular, if $a=1$ then our construction coincides with the lexicographic
product of $\Gamma$ with the complete graph $K_{j+1}$. 
The similarity with X-joins can be seen after collapsing each cell in $\Gamma^j_{\pi}$
into a single vertex. From this follows that the only difference between our construction
and X-join is, that we define differently adjacencies between vertices in the cells 
which do belong to the same copy of $\Gamma$, and those which do not.

Now, we translate the $\pi$-join construction into the language of matrices, 
since it will be easier to prove statements using adjacency matrices of graphs.

Let $\Gamma$ be a (di)graph and $\pi=\{C_1,C_2,\ldots,C_a\}$ be a homogeneous 
partition of its vertices into $a$ cells of size $b$, for suitable positive 
integers $a$ and $b$. Reorder the vertices of $\Gamma$ according to $\pi$ 
and create adjacency matrix $A$ of $\Gamma$ which respects $\pi$. 
Then we may consider $A$ as a block matrix with blocks of size $b\times b$.
Let us define for $i\in\{1,\ldots,a\}$ matrix $U_i$ as 
a matrix whose all entries are ones in the columns with indices
$(i-1)a+1, (i-1)a+2, \ldots, i\cdot a$, and all other entries
are zeros, i.e. $U_i=(0,\ldots,0,1,0,\ldots,0)\otimes J$,
where $J$ is the all-ones matrix of size $n\times b$, 
and $\otimes$ stands for the Kronecker product. 
Further, for an arbitrary integer $j$ let us create 
a circulant block matrix $M_j(A)$ whose first row of blocks is 
$$(A, \underbrace{U_1,\ldots, U_1,}_{j} \underbrace{U_2,\ldots,U_2}_{j}, \ldots, 
\underbrace{U_a,\ldots, U_a}_{j}),$$ 
and all other rows of blocks are obtained by shifting the blocks to right by one position.   

The matrix $M_j(A)$ is the adjacency matrix of the $\pi$-join of $\Gamma$ in 
power $j$, defined above.

\begin{example}
For $n=6$, $a=3$ and $b=2$ then we have the following matrices:
$$
U_1=\left(
\begin{array}{cc|cc|cc}
1 & 1 & 0 & 0 & 0 & 0 \\
1 & 1 & 0 & 0 & 0 & 0 \\
1 & 1 & 0 & 0 & 0 & 0 \\
1 & 1 & 0 & 0 & 0 & 0 \\
1 & 1 & 0 & 0 & 0 & 0 \\
1 & 1 & 0 & 0 & 0 & 0 \\
\end{array}
\right),
\,
U_2=\left(
\begin{array}{cc|cc|cc}
0 & 0 & 1 & 1 & 0 & 0 \\
0 & 0 & 1 & 1 & 0 & 0 \\
0 & 0 & 1 & 1 & 0 & 0 \\
0 & 0 & 1 & 1 & 0 & 0 \\
0 & 0 & 1 & 1 & 0 & 0 \\
0 & 0 & 1 & 1 & 0 & 0 \\
\end{array}
\right),
\,
U_3=\left(
\begin{array}{cc|cc|cc}
0 & 0 & 0 & 0 & 1 & 1 \\
0 & 0 & 0 & 0 & 1 & 1 \\
0 & 0 & 0 & 0 & 1 & 1 \\
0 & 0 & 0 & 0 & 1 & 1 \\
0 & 0 & 0 & 0 & 1 & 1 \\
0 & 0 & 0 & 0 & 1 & 1 \\
\end{array}
\right),
$$
while for $j=1$ and $j=2$ we have:
$$
M_1(A)=\left(
\begin{array}{cccc}
A & U_1 & U_2 & U_3 \\
U_3 & A & U_1 & U_2 \\
U_2 & U_3 & A & U_1 \\
U_1 & U_2 & U_3 & A \\
\end{array}
\right),
\qquad
M_2(A)=\left(
\begin{array}{ccccccc}
A & U_1 & U_1 & U_2 & U_2 & U_3 & U_3\\
U_3 & A & U_1 & U_1 & U_2 & U_2 & U_3\\
U_3 & U_3 & A & U_1 & U_1 & U_2 & U_2\\
U_2 & U_3 & U_3 & A & U_1 & U_1 & U_2\\
U_2 & U_2 & U_3 & U_3 & A & U_1 & U_1\\
U_1 & U_2 & U_2 & U_3 & U_3 & A & U_1\\
U_1 & U_1 & U_2 & U_2 & U_3 & U_3 & A\\
\end{array}
\right).
$$
\end{example}

\subsection{Motivation}

The motivation behind the $\pi$-join construction defined in the previous subsection 
comes from an analysis of previously known DSRGs. 
It showed that a significant amount of them 
can be constructed as a $\pi$-join of certain smaller DSRG $\Gamma$, 
its partition $\pi$ and a suitable integer~$j$. 
Inspired by this observation we started to investigate the following problem:

\medskip

\noindent{\bf Problem:} {\it 
For a given (D)SRG $\Gamma$, what are the necessary and sufficient conditions
for a homogeneous partition $\pi$ in order to get a DSRG from the $j$-th power of its $\pi$-join?
}

\medskip

We will solve this problem in two steps. First, we will consider the smallest case $j=1$ 
to derive some necessary conditions on $\Gamma$ and $\pi$. After that 
we will show that if these conditions are satisfied, then for arbitrary integer
$j$ the $j$-th power of the $\pi$-join produces again a DSRG. If a $\pi$-join
of a (D)SRG $\Gamma$ produces a larger DSRG, then we say that $\pi$ is
a \emph{good partition} for $\Gamma$.  

\subsection{Necessary conditions}

Here we derive the necessary conditions for a homogeneous partition $\pi$ in 
a (D)SRG $\Gamma$ in order to be good for the $\pi$-join construction.
Let us start with an easy observation about matrices $U_i$.
Computational arguments lead to the following properties:

\medskip

\noindent{\bf Observation 1: } Let $a,b,k,A,U_i$ are as in the definition of the 
$\pi$-join construction. Then for any $i\in\{1,\ldots,a\}$ and $l\in\{1,\ldots,b\}$ we
have:
\begin{itemize}
\item[(i)] $\sum_{i=1}^aU_i=J$,
\item[(ii)] $U_i\cdot U_l=b\cdot U_l$,
\item[(iii)] $A\cdot U_l=k\cdot U_l$,
\item[(iv)] $U_i\cdot A$ has constant columns, and this constant in an
arbitrary column represents the number of darts starting from cell $C_i$ 
and terminating in the vertex corresponding to the given column. 
\end{itemize}

The following theorem is very important, since it has many consequences, which
significantly restricts our search space for possibly good partitions.

\begin{theorem}
Let $A$ be an adjacency matrix of a (D)SRG $\Gamma$ with parameters $(n,k,t,\lambda,\mu)$, 
which respects the homogeneous partition $\pi=\{C_1,\ldots,C_a\}$ of degree $b$.
Suppose that the $\pi$-join $\Gamma^1_{\pi}$ for $j=1$ is a DSRG
with parameters $(\tilde n,\tilde k,\tilde t,\tilde\lambda,\tilde\mu)$.
Then 
\begin{itemize}
\item[(a)] $\tilde n=(a+1)n,\quad \tilde k=n+k,\quad \tilde t=b+t,\quad \tilde\lambda=b+\lambda,\quad 
\textrm{ and }\quad\tilde\mu=b+\mu.$
\item[(b)] For arbitrary $i,l\in\{1,\ldots,a\}$ the number $q_{i,l}$ of darts
starting in $C_i$ and terminating in $v\in C_l$ does not depend
on the concrete choice of $v$, just on $i$ and $l$. Moreover,
$$q_{i,l}=\begin{cases} \lambda+b-k &\textrm{ if } i=l, \\ \mu &\textrm{ if } i\neq l.\end{cases}$$
\end{itemize}
\end{theorem}

\noindent{\bf Proof.} 
Denote $M_1(A)$ the adjacency matrix of $\Gamma^1_{\pi}$ coming from the construction.
Clearly, $M_1^2(A)$ will be a circulant block matrix which is completely determined
by its first row of blocks, so let us say that its first row is $(B_0, B_1, \ldots, B_a)$. 

\begin{itemize}
\item[(a)]
The equations for $\tilde n$ and $\tilde k$ follow immediately from the method
of construction. For other parameters let us compute:
$$ B_0 = A^2 + U_1U_{a} + U_2U_{a-1} + \ldots + U_aU_1 =
A^2+b\cdot (U_a+ U_{a-1}+\ldots + U_1) = $$
$$= A^2 + b\cdot J = \lambda\cdot A+\mu(J-I-A)+t\cdot I + b\cdot J=
(\lambda+b)A+(\mu+b)(J-I-A)+(t+b)I,$$
and the claim follows.
\item[(b)] Again compute:
$$ B_1 = AU_1+U_1A+U_2U_a+\ldots+U_aU_2= U_1A+kU_1+b(J-U_1), $$
but from strong regularity we have also
$ B_1 = \tilde\lambda U_1+\tilde\mu(J-U_1). $
Altogether this means that $U_1A=(\lambda+b-k)U_1+\mu(J-U_1)$. 
Now using part $(iv)$ from {\it Observation 1} we get, that
the number of darts starting in $C_1$ and terminating in a concrete 
vertex in $C_1$ is precisely $\lambda+b-k$; while the number of 
those starting from $C_1$ and terminating in a concrete vertex in $C_i$, 
where $i>1$, is precisely $\mu$. Analogous computations with $B_2,B_3,
\ldots,B_a$ prove the claim.   
\end{itemize}

\rightline{$\square$}

\medskip

\noindent{\bf Corollary 1.}
Let $M_1(A)$ and $A$ are as above. Then if $M_1(A)$ is an adjacency matrix 
of a DSRG then necessarily 
\begin{equation}
2k+\mu-\lambda=a\mu+b.\tag{eq1}
\end{equation}

\medskip

\noindent{\bf Proof.} 
Counting the number of incoming darts to a fixed vertex in $C_l$ in $\Gamma$ 
we get $$k=\sum_{i=1}^aq_{i,l}=(\lambda+b-k)+(a-1)\mu.$$

\medskip

\noindent{\bf Corollary 2.}
Necessarily:
\begin{equation}
\lambda+b-k\geq 0. \tag{eq2}
\end{equation}

\medskip

\noindent{\bf Corollary 3.} 
If $M_1(A)$ is an adjacency matrix of a DSRG, then $\pi$ is a homogeneous 
column-equitable partition with quotient matrix $Q=(\lambda+b-k)I+\mu (J-I)$.

\medskip

\noindent{\bf Proof.} 
Immediately follows from Corollary 1.

\medskip

\noindent{\bf Corollary 4.}
If $A$ is an adjacency matrix with respect to a partition $\pi$ of an undirected
strongly regular graph, then necessarily $\pi$ is an equitable partition
and the cells are regular graphs of valency $\lambda+b-k$.

\medskip

\noindent{\bf Proof.} In an undirected graph indegree and outdegree of a vertex are equal, 
therefore the partition $\pi$ satisfies all conditions for being equitable. 

\medskip

\noindent{\bf Corollary 5.}
For given (D)SRG with parameter set $(n,k,t,\lambda,\mu)$ the number of 
solutions of the equation \thetag{eq1} is at most two. 

\medskip

\noindent{\bf Proof.} 
Since $ab=n$, replacing $a$ by $n/b$ results in a quadratic equation in $b$.

\medskip

\noindent{\bf Corollary 6.}
If $a=1$, then $b=n$, and the potential parameter set of the 
$\pi$-join is $(2n,n+k,n+t,n+\lambda,n+\mu)$. Feasibility conditions for parameter sets
imply that $n=1$, which is not an interesting case. 

\medskip

\noindent{\bf Proof.} 
If $(2n,n+k,n+t,n+\lambda,n+\mu)$ is a parameter set of a DSRG, then
from (\ref{c1}) we have
\[
(n+k)(n+k+n+\mu-n-\lambda)=(t+n)+(2n-1)(\mu+n).
\]
After simple algebraic manipulations, using equations 
$k(k+\mu-\lambda)=t+(n-1)\mu$ and \thetag{eq1} we obtain that $n=1$. 

\medskip

\noindent{\bf Corollary 7.} 
If $b=1$, then the initial graph $\Gamma$ is necessarily complete.

\medskip

\noindent{\bf Proof.} The cells are single vertices, therefore part (b) of Lemma~1 
implies that the initial graph is necessarily complete. 

\subsection{Sufficient conditions}

In this section we state that the conditions derived in the previous 
subsection are, in fact, sufficient. 

\begin{theorem}
Let $\Gamma$ be a (D)SRG with parameter set $(n,k,t,\lambda,\mu)$. Let $a$ and $b$ are 
positive integers such that $ab=n$ and there exists a homogeneous column 
equitable partition $\pi=\{C_1,\ldots,C_a\}$ 
of vertices with quotient matrix $Q=(\lambda+b-k)I+\mu(J-I)$. 
Let $A$ be an adjacency matrix respecting $\pi$, and let us 
define matrix $M_j(A)$ in accordance with our $\pi$-join construction for 
an arbitrary positive integer $j$. 
Then $M_j(A)$ is an adjacency matrix of a DSRG
with parameter set $$( (ja+1)n, jn +k, jb+t, jb+\lambda, jb+\mu ).$$ 
\end{theorem}

\noindent{\bf Proof.} The necessity of conditions has been shown above. By direct computations 
one can show that if $A$ satisfies all the posed conditions, 
then they are also sufficient, therefore $M_j(A)$ 
corresponds to a DSRG for any positive integer~$j$.

\medskip

This theorem provides an algorithm how to proceed when we are looking for 
a $\pi$-join of a DSRG:
\begin{enumerate}
\item Take a directed or an undirected SRG $\Gamma$; 
\item Solve \thetag{eq1} for the parameter set of $\Gamma$;
\item For the solutions of \thetag{eq1} satisfying \thetag{eq2} count the 
quotient matrix and find a column equitable partition $\pi$ 
with this quotient matrix. 
\item For $\pi$, $\Gamma$ and an arbitrary integer $j$ create the 
$\pi$-join in power $j$ for $\Gamma$ according to the construction.  
\end{enumerate}

Of course, the key moment of the construction is the finding of a good 
column equitable partition. In the last sections we try to provide such partitions
(if they do exist) at least for the small cases of DSRGs and for some general
families. As we will see, in some cases the method of construction
of the initial DSRG immediately provides us such a partition. 
In other cases some combinatorial tricks help to find it; 
and it happens sometimes that we are forced to leave this job to a computer.  

\subsection{Properties of the construction}

Here we investigate two properties of the $\pi$-joins: 
behaviour under the complementation, and the set of eigenvalues. 

\subsubsection{Complementation}

In this subsection we investigate when does there exist a 
partition $\pi$ for which both (D)SRG and its complementary (di)graph 
can be $\pi$-joined into a new DSRG, simultaneously.   
After that we derive the sufficient condition when we can get 
two DSRGs with complementary parameter sets as a result of the 
$\pi$-join construction.

\begin{theorem}
Suppose that $\Gamma$ is a (D)SRG with parameters $(n,k,t,\lambda,\mu)$, 
$\bar \Gamma$ its complementary digraph with parameters 
$(n,\bar k, \bar t,\bar\lambda,\bar\mu)$
and $\pi$ a homogeneous column equitable partition for both $\Gamma$ and $\bar\Gamma$.
Then $\mu+\bar\mu=n/2$ and $\pi$ has precisely two cells. 
\label{thm2} 
\end{theorem}

\medskip

\noindent{\bf Proof.}
From the necessary conditions on the quotient matrix follow that 
$\tilde\mu+\mu=b$ and $(\lambda+b-k)+(\tilde\lambda+b-\tilde k)=b-1$. 
Using the identities for the parameters of the complementary digraph of a DSRG
we get: $n-2k+\lambda+\mu=b$ from the first equation. After this the second
equation gives us:
$$b-1=(\lambda+\tilde\lambda)+2b-(k+\tilde k)=\lambda+n-2k+\mu-2+2b-(n-1)=$$
$$=b-2+2b-n+1=3b-1-n,$$
hence $n=2b$ and the claim follows.

\medskip

From the proof and the assumptions of the previous theorem we can derive the 
following corollary:

\medskip

\noindent{\bf Corollary 8.} 
Let $\Gamma$ be a DSRG with parameters $(n,k,t,\lambda,\mu)$,
$\bar\Gamma$ its complement with parameters $(n,\bar k, \bar t,\bar\lambda,\bar\mu)$.
Suppose that $\mu+\tilde\mu=n/2$ and $\pi$ is a good homogeneous column-equitable 
partition with two cells for $\Gamma$. 
Then it is a good partition for $\bar\Gamma$ as well.

\medskip

\begin{theorem}
Let $\Gamma$ be a DSRG with parameters $(n,k,t,\lambda,\mu)$,
$\bar\Gamma$ its complement with parameters $(n,\bar k, \bar t,\bar\lambda,\bar\mu)$.
Suppose that they can be simultaneously $\pi$-joined via partition $\pi$ to a DSRG.
Then $\Gamma^{\pi}$ and $\bar\Gamma^{\pi}$ have complementary parameter sets. 
\label{thm3}
\end{theorem}

\medskip

\noindent{\bf Proof.}
We already know from Theorem~\ref{thm2} that under the posed assumptions the partition
$\pi$ has two cells. Now putting $a=2$ into the formulas for parameter sets
of the $\pi$-join and using relations between parameter sets of mutually complementary
digraphs we get that the parameter sets for $\overline{\Gamma^{\pi}}$ and 
$\bar\Gamma^{\pi}$ coincide. 

\subsubsection{Eigenvalues}

It is well-known that a DSRG has three different eigenvalues
and all of them are integers. These eigenvalues are uniquely determined
by the parameter set. 
It is not hard to see that if the basic graph $\Gamma$ with
parameter set $(n,k,t,\lambda,\mu)$ has eigenvalues 
$\{k,\theta_1,\theta_2\}$ then its $\pi$-join 
has eigenvalues $\{jn+k, \theta_1, \theta_2\}$. 

\section{Good partitions for general families of (D)SRGs}

In this section we try to describe some good equitable partitions
for well-known families of (D)SRGs. Consequently, the $\pi$-join construction
gives us many infinite families of DSRGs, which are in many cases new.
For fixed parameter set we consider just one graph and we focus on finding 
at least one good equitable partition. 

\subsection{Complete graphs}
\label{ss7.1}

The complete graph $K_n$ can be regarded as an SRG with parameters $(n,n-1,n-2,\mu)$, 
where the parameter $\mu$ is undefined. Thus, we cannot apply our previous theory. 
However, if we apply the $\pi$-join construction for the complete graph $K_{ab}$ 
using the homogeneous equitable partition with $a$ cells of size $b$, 
then it is easy to show that the result is again a DSRG. Its parameter set is 
$$
(ja^2b+ab, jab+ab-1, jb+ab-1, jb+ab-2, jb+b).
$$
This parameter set is known for example from \cite{os}, Section 4.

\subsection{Duval's construction for $\lambda=0$, $t=\mu$}
\label{ss7.2}

Duval (in \cite{du}, Lemma 8.2) constructed an infinite family of DSRGs 
with parameter set $(k^2+k,k,1,0,1)$ by giving its adjacency matrix. 

It is an easy exercise to check that his matrix defines a good partition
for both feasible integer solutions of \thetag{eq1}: $(a,b)=(k,k+1)$ and
$(a,b)=(k+1,k)$. More precisely, if the vertices $v_1,v_2, \ldots, v_n$
of the graph are in accordance with the order of rows/columns of the 
adjacency matrix given by Duval, then 
$$\pi_1=\left\{\{v_1,v_2,\ldots,v_{k+1}\},\{v_{k+2},\ldots,v_{2k+2}\},
\ldots,\{v_{k^2},v_{k^2+1},\ldots,v_{k^2+k}\}\right\},$$ and 
$$\pi_2=\left\{\{v_1,v_2,\ldots,v_k\},\{v_{k+1},\ldots,v_{2k}\}, \ldots,
\{v_{k^2+1},\ldots,v_{k^2+k}\right\}\}$$
are the corresponding good partitions. 

This proves the existence of DSRGs with the following parameter sets:
\[
\left(k(k+1)(jk+1), jk(k+1)+k, jk+j+1, j(k+1), jk+j+1\right),
\]
and
\[
\left(k(k+1)(jk+j+1), jk(k+1)+k, jk+1, jk, jk+1\right),
\]
where $j$ and $k$ are arbitrary positive integers. 
Among these parameter sets there are several new: 
(60,15,4,3,4), (60,20,7,6,7), (78,26,9,8,9), (96,32,11,10,11), and (108,27,7,6,7)
for $n\leq 110$. 

In similar manner one can show that from Duval's generalized construction 
(\cite{du}, Theorem 8.3),
in conjunction with our construction it is possible to obtain DSRGs with parameter sets

\[
\left(mk(k+1)(jk+1), mjk(k+1)+mk, m(jk+j+1), mj(k+1), m(jk+j+1)\right),
\]
\[
\left(mk(k+1)(jk+j+1), mjk(k+1)+mk, m(jk+1), mjk, m(jk+1)\right),
\]
for arbitrary positive integers $j,m,k$. 

\subsection{J\o rgensen's construction}
\label{ss7.3}

For any positive integers $k$ and $\mu$ such that $\mu\mid k-1$, J\o rgensen 
in \cite{jg} constructed a DSRG $\Gamma$ with parameter set $((k^2-1)/\mu, k, \mu+1,
\mu,\mu)$. The vertices of his graphs are integers modulo $n=(k^2-1)/\mu$,
and there is a dart from $x$ to $y$ if and only if $x+ky\in\{1,2,\ldots,k\}$.
(Operations are also taken modulo $n$.) In other words, a fixed vertex $x$ has 
outneighbours of form $k(s-x)$ and inneighbours of form $s-ky$ for $s=1,2,\ldots,k$. 

For this parameter set \thetag{eq1} has always a solution in integers:
$(a,b)=((k-1)/\mu, k+1)$. 
The possibly second solution is $(a,b)=((k+1)/\mu,k-1)$, but since $\mu\mid k-1$
an easy number theoretical argument shows that $(k+1)/\mu$ is integer if and 
only if $\mu\in\{1,2\}$. 
  
\subsubsection{Case $a=\frac{k-1}{\mu}$, $b=k+1$}

Let us define for $i\in\{0,1,\ldots,a-1\}$ the following sets: 
$$C_i=\{bi\}\cup \{k(s-bi)\,|\, s\in\{1,2,\ldots,k\}\},$$
i.e. $C_i$ contains the number $bi$ and its outneighbours.
It is not hard to show that if from vertex $bi$ there is a dart
to some $y$, then there is no dart from $bj$ to $y$ for $j\neq i$,
and there is no dart between $bi$ and $bj$. Therefore the sets $C_i$
are cells of a partition of $\{0,1,2,\ldots,n-1\}$
which was identified with the vertex set of $\Gamma$. 

Pick a vertex $r\neq b\cdot i$ in $C_i$. Then, clearly, $bi\to r$ is a dart in $\Gamma$,
and we know that there are precisely $\lambda$ directed two-paths from $bi$ to $r$
in $\Gamma$. However, all these two-paths are inside the cell $C_i$, because
all the outneighbours of $bi$ belong to $C_i$. Hence, the number of inneighbours
of $r$ starting from $C_i$ is precisely $\lambda+1$, and the number of inneighbours 
of $bi$ in $C_i$ is also $\lambda+1$, which follows from the property that $\Gamma$ is
a DSRG with parameter $t=\lambda+1$. Altogether this means that for each vertex in $C_i$
there are precisely $\lambda+1$ inneighbours in $C_i$, thus if our partition is column equitable,
then $q_{i,i}=\lambda+1=\mu+1$. 

Now take $i\neq j$, an arbitrary vertex $r\in C_j$, and count the number of inneighbours
of $r$ in $C_i$. Clearly, $bi$ can not be an inneighbour of $r$, since $r\notin C_i$.
Therefore there are precisely $\mu$ oriented paths of length two from $bi$ to $r$, 
and we know that the midvertices on these paths lie in $C_i$. Hence, $q_{i,j}=\mu$. 

In fact, we proved the following lemma:

\begin{lemma}
For arbitrary $k$ and $\mu$ satisfying $\mu\mid k-1$
let $a=(k-1)/\mu$, $b=k+1$ and $\Gamma$ be the corresponding
DSRG from J\o rgensen's construction.  
For $i\in\{0,1,\ldots,a-1\}$ the sets 
$$C_i=\{ib\}\cup \{k(s-ib)\,|\, s\in\{1,2,\ldots,k\}\}$$
form a column equitable partition of $\Gamma$ with quotient matrix
$Q=(\mu+1)I+\mu(J-I)$.
Therefore for an arbitrary positive integer $j$ there exists a 
DSRG with parameter set $(N,K,\Lambda+1,\Lambda,\Lambda)$, where:
$$
N=\left(\frac{j(k-1)}{\mu}+1\right)\cdot\frac{k^2-1}{\mu},\quad
K=\frac{j(k^2-1)}{\mu}+k,\quad
\Lambda=j(k+1)+\mu.
$$
\end{lemma}  

Calculations show that $\Lambda\mid K-1$ does hold, therefore this parameter set
is not new, since it is covered by the construction given by J\o rgensen. However, 
the resulting graphs from different constructions are not isomorphic, in general.
For example, when $j=1$, $k=7$ and $\mu=3$ then $K=23$ and $\Lambda=11$, 
but the resulting DSRGs from our and J\o rgensen's
construction with parameters $(48,23,12,11,11)$ are non-equivalent.  

\subsubsection{Case $a=\frac{k+1}{\mu}$, $b=k-1$}

Recall that in this case the assumption $\mu\mid k-1$ on the initial parameters 
imply that the only possibilities are $\mu\in\{1,2\}$, otherwise $a$ cannot be an
integer. 
Now let us consider a partition of the set $\{0,1,2,\ldots, n-1\}$ 
into cells $C_i$ as equivalence classes modulo $a$, i.e.
$$C_i=\{s \quad|\quad s\equiv i \pmod a\}.$$
We show that $\{C_0,C_1,\ldots,C_{a-1}\}$ is a good equitable partition. 

Let us select any $y\in C_i$ and count the number of $x\in C_i$ such that
$x\to y$ is a dart. We repeat that if $x$ is an inneighbour of $y$, 
then for some $s\in\{1,2,\ldots,k\}$ we have $x=s-ky$, 
therefore necessarily $s-ky\equiv y \pmod a$. Equivalently,
$s\equiv (k+1)y\equiv 0\pmod a$, because $k+1$ is a multiple of $a$. 

\begin{itemize}
\item If $\mu=1$, then $s\equiv 0\pmod{k+1}$ has
no solution in $s\in\{1,2,\ldots,k\}$. 
\item If $\mu=2$, then $s\equiv 0\pmod{(k+1)/2}$ has a unique solution $s=(k+1)/2$.  
\end{itemize}

Now let us consider the number of inneighbours in $C_j$ of a fixed vertex $y\in C_i$,
where $i\neq j$. Similarly as above, for an inneighbour $x\in C_j$ there is $s\in\{1,2,\ldots,k\}$
such that $x=s-ky$, and $y-x\equiv j-i\pmod a$. After simplification: $s\equiv i-j\pmod a$. 

\begin{itemize}
\item If $\mu=1$, then $s\equiv i-j\pmod{k+1}$ has a unique solution $s=i-j$ (taken modulo $a$). 
\item If $\mu=2$, then $s\equiv i-j\pmod{(k+1)/2}$ has exactly two solutions.
\end{itemize}

We proved the following:

\begin{lemma}
For $\mu\in\{1,2\}$ and arbitrary $k$ satisfying $\mu\mid k-1$
let $a=(k+1)/\mu$, $b=k-1$ and $\Gamma$ be the corresponding
DSRG from J\o rgensen's construction.  
For $i\in\{0,1,\ldots,a-1\}$ the sets: 
$$C_i=\{s \quad|\quad s\equiv i \pmod a\}$$
form a column equitable partition of $\Gamma$ with quotient matrix
$(\mu-1)I+\mu(J-I)$.
Therefore for arbitrary positive integer $j$ there exists a 
DSRG with parameter set $(N,K,\Lambda+1,\Lambda,\Lambda)$, where:
$$
N=\left(\frac{j(k+1)}{\mu}+1\right)\cdot\frac{k^2-1}{\mu},\quad
K=\frac{j(k^2-1)}{\mu}+k,\quad
\Lambda=j(k-1)+\mu.
$$
\end{lemma}

In the case, when $k=1$ the resulting parameter set can be also
obtained from J\o rgensen's construction. Also this is the case, when $k=2$, $\mu=j=1$.
In all other cases $\Lambda$ does not divide $K-1$. According to
\cite{ab} we constructed first DSRG with parameters $(40,11,4,3,3)$, $(72,19,6,5,5)$,
$(90,19,5,4,4)$ and $(104,27,8,7,7)$, if we consider just $N\leq 110$. 

\subsection{Triangular graphs}
\label{ss7.4}

The \emph{triangular graph $T(n)$} is defined as the line graph of the complete
graph $K_n$. For $n\geq 4$ it is an SRG with parameters $({n\choose 2}, 2(n-2), n-2, 4)$.
Moreover, graphs with these parameter set are unique, when $n\neq 8$. 
For $n=8$ there are three more graphs, the so-called \emph{Chang graphs} 
(see \cite{cc, ch}), apart from $T(8)$. 

By solving $\thetag{eq1}$ we get that for $n\equiv 1\pmod 2$ there is 
a unique solution $(a,b)=(\frac{n-1}2,n)$; for $n\equiv 0\pmod 4$ also a unique solution
$(a,b)=(n/4, 2n-2)$, while for $n\equiv 2\pmod 4$ there is no integer solution.

\subsubsection{Case $n\equiv 1\pmod 2$}
For odd $n$ we have $a=\frac{n-1}2$ and $b=n$. So the cells of the potential equitable
partitions are regular of valency $\lambda+b-k=2$. In other words, each cell is 
a collection of cycles. 

Let us consider the vertex set of the triangular graph as 2-subsets of the set
$\mathbb Z_n=\{0,1,2,\ldots,n-1\}$ and there is an edge between two different 
2-subsets if and only if they are not disjoint.
For $i=1,2,3,\ldots, (n-1)/2$ let us define sets 
$C_i=\{\{\ell,\ell+i\}\,|\, \ell\in\mathbb Z_n\}$. 
The partition $\pi=\{C_1,C_2,\ldots,C_{(n-1)/2}\}$ is equitable with quotient matrix
$Q=2I+4(J-I)$, thus the existence of DSRGs with parameter set
$(N,K,T,\Lambda,M)$ is proved, where 
$N=m(2m+1)(jm+1)$, $K=jm(2m+1)+4m-2$, $T=j(2m+1)+4m-2$, $\Lambda=j(2m+1)+2m-1$, 
$M=j(2m+1)+4$
and $j,m$ are arbitrary integers. 
In all the cases when $N\leq 110$ the parameter set has been known previously.

\subsubsection{Case $n\equiv 0\pmod 4$}
Here the cells are of order $2n-2$ and valency $\lambda+b-k=n$. Again, equivalently we can
consider the graph $K_n$ and correspondingly an edge partition in it. The partition determines
an edge-coloring (not proper) in the sense, that the same color is assigned to edges from the 
same cell and different colors to edges from different cells. Each color defines a 
subgraph of $K_n$. We can assume that the edge $u_1u_2$ is black and we further consider the 
``black subgraph'' of $K_n$. This edge is incident to $n$ black edges, therefore 
$\deg_B(u_1)+\deg_B(u_2)=n+2$, where $\deg_B(u)$ stands for ``black-degree'', i.e. the number 
of black edges incident to $u$. By pigeonhole principle there exists a vertex $u_3$ such that
$u_1u_3$ and $u_2u_3$ are black edges. Thus, 
$\deg_B(u_1)+\deg_B(u_3)=\deg_B(u_2)+\deg_B(u_3)=n+2$, hence 
$\deg_B(u_1)=\deg_B(u_2)=\deg_B(u_3)=1+\frac{n}2$, and necessarily all the $n/2+1$ neighbours of
$u_1$ have black-degree $n/2+1$. So the number of black edges is at least 
$\frac12(1+\frac{n}2)(2+\frac{n}2)$. However, for $n>6$ this number is greater than $2n-2$, 
thus for $n\geq 8$ there is no good partition. For $n=4$ the partition is trivial. 

For the three exceptional cases, Chang graphs, we executed an exhaustive computer search
for finding induced regular subgraphs of size 14 and valency 8. Since the number of subsets
of size 14 in a set of size 28 is quite large, we used the {\sf SetOrbit} package written
by Pech and Reichard \cite{sr}. 
With its aid we can generate just one representative for each orbit of the
group of automorphisms on subsets of size 14. 
This tool reduced the amount of subsets being considered for a few hundreds of thousands, 
and therefore it took just a few hours of computation. 
However, this exhaustive search showed that there is no such induced 
subgraph, and therefore there is no good partition for our construction.   

\subsection{Complements to triangular graphs}
\label{ss7.5}

Triangular graph $T(n)$ has been described above. Let us now consider its complement $\bar T(n)$. 
Again, to avoid trivial cases we suppose that $n>4$.  The parameter set of $\bar T(n)$ is
$\left( {n\choose 2}, {n-2\choose 2}, {n-4\choose 2}, {n-3\choose 2}\right)$. Equation \thetag{eq1}
has integer solution just for $n=5,6,8$. For $n=5$ see the section about the Petersen graph. For
$n=6$ exhaustive computer search showed that there is no good partition for the only solution 
$(a,b)=(3,5)$. Finally, for $n=8$ we have just one solution $(a,b)=(2,14)$, but $\mu+\bar\mu=14$ 
and Theorem~\ref{thm2} applies.

\subsection{Lattice square graphs}
\label{ss7.6}

Another famous family of strongly regular graphs is the so-called \emph{lattice square graphs}. 
The lattice square graph $L_2(n)$ has vertex set $V=\{1,2,\ldots,n\}\times\{1,2,\ldots,n\}$
and two vertices are adjacent if their coordinates agree in exactly one position.
The parameter set of $L_2(n)$ is $(n^2, 2n-2, n-2, 2)$. 
For each positive integer $n$ the equation \thetag{eq1} has at least one solution $a=b=n$. 
For odd $n$ this is the only solution, but for even $n$ there is a second 
solution $(a,b)=(n/2,2n)$.  

Conditions posed on the quotient matrix $Q$ of a good partition say that in the first
case $Q=2(J-I)$, and in the second $Q=nI+2(J-I)$.

First, consider partition $\pi=\{C_1,C_2,\ldots,C_n\}$ of vertices of $L_2(n)$ such that

\[
\begin{array}{rcl}
C_1 & = & \{ (1,1), (2,2), \ldots, (n-1,n-1), (n,n)\}, \\
C_2 & = & \{ (1,2), (2,3), \ldots, (n-1,n), (n,1)\}, \\
& \vdots & \\
C_n & = & \{ (1,n), (2,1), \ldots, (n-1,n-2), (n,n-1)\}.\\
\end{array}
\]

It is not hard to see that this partition is equitable with quotient matrix $Q=2(J-I)$,
and it guarantees the existence of a DSRG with parameters
\[
\left(n^2(jn+1), jn^2+2n-2, jn+2n-2, jn+n-2, jn+2\right)
\]
for arbitrary positive integers $j$ and $n$. The resulting parameter set 
is new at least in the case when $j=1$ and $n=4$. 

When $n$ is even, we can consider the partition 
$\pi=\{C_1,C_2,\ldots,C_{n/2}\}$ such that for all $m=1,2,\ldots,n/2$
the cell $C_m$ is: 
\[
C_m=\left\{(x,y)\,|\, x\in\{1,2,\ldots,n\}, y\in\{2m-1,2m\}\right\}.
\]
It is an equitable partition with quotient matrix $Q=nI+2(J-I)$. 
The corresponding $\pi$-join of the lattice square graph $L_2(n)$
is a DSRG for each integer $j$ with parameter set 
\[
\left(n^2(1+\frac{jn}2), jn^2+2n-2, 2jn+2n-2, 2jn+n-2, 2jn+2\right).
\]

\subsection{Construction by Klin et al.}

Klin et al. in \cite{km} constructed a family of DSRGs with parameter set 
$(2n,n-1,\frac{n-1}2, \frac{n-3}2, \frac{n-1}2)$, where $n$ is odd, 
as Cayley graphs over dihedral group of order $2n$.
So let $D_n=\langle r,s\,\mid\, r^n=s^2=1\rangle$ be the dihedral group
of order $2n$. Let $X=\{r,r^2,\ldots,r^{(n-1)/2},sr, sr^2,\ldots, sr^{(n-1)/2}\}$.
The Cayley graph with connection set $X$ over $D_n$ is the requested DSRG.

For this parameter set \thetag{eq1} has always a unique solution, when $n>3$. 
This is $(a,b)=(2,n)$. For a good partition the quotient matrix is
$(n-1)/2\cdot J$, and it is not hard to see that the partition
$$\pi=\{\{1,r,r^2,\ldots,r^{n-1}\},\{s,sr,sr^2,\ldots,sr^{n-1}\}\}$$
satisfies this condition, therefore $\pi$ is a good partition. 
Hence the existence of DSRGs with parameter set 
$$\left((4j+2)n, 2jn+n-1, jn+\frac{n-1}2, jn+\frac{n-3}2, jn+\frac{n-1}2\right)$$ 
is proved for arbitrary $j$ and odd $n$. 

\subsection{Construction by Duval and Iourinski}
\label{ss7.7}

Duval and Iourinski in \cite{di} gave a sufficient condition for construction DSRGs as Cayley graphs
over metacirculant groups. For our construction it is sufficient to choose their orbit partition
satisfying the $q$-orbit condition in order to get a good partition for $(a,b)=(q,m)$.

For their graphs equation \thetag{eq1} has a second solution just in the case, when $m=q+1$. 
This case coincides with the construction leading to DSRGs with $\mu=0$, $t=\lambda=1$.

\section{Good partitions for small cases of (D)SRGs} 

In this section we discuss the $\pi$-joins of small 
(D)SRGs with parameter sets which have not been included in the previous section. 
We show just a single good partition for each parameter set, usually. 
In fact, in the smallest cases ($n\leq 16$) we run an exhaustive search with the
aid of a computer for the values $j=1,2,3$. In many cases we have found 
dozens of good partitions and non-isomorphic $\pi$-joins. 
These are available upon request from the author. 

Exceptionally, we will consider more graphs sharing the same parameter set or a 
graph of larger order. It happens since we want to include into 
consideration some famous graphs, as well. For example, the Petersen graph, 
the Shrikhande graph, the Schl\"afli graph or the Hoffman-Singleton graph.
On the other hand, we will not consider parameter sets $(n,k,t,\lambda,\mu)$,
where $k\geq n/2$ and the only solution of the \thetag{eq1} is the pair $(a,b)=(2,n/2)$,
since in these cases we can apply Theorem~\ref{thm2} and Theorem~\ref{thm3}.

A summary of the all new parameter sets of DSRGs is given in Appendix.

\subsection{SRG(4,2,0,2)}

The smallest undirected SRG is a 4-cycle. The only good solution of \thetag{eq1}
is $a=b=2$. Though there are two equitable partitions into two cells of size 2,
just one of them, the pairs of non-adjacent vertices, satisfies the condition posed 
on the quotient matrix. 
Our construction creates DSRGs with parameter set $(8j+4,4j+2,2j+2,2j,2j+2)$.
For this parameter set and its complementary parameter set there are known 
several constructions, see \cite{du}, \cite{hs}, \cite{jg} and \cite{os}. 

\subsection{SRG(10,3,0,1) -- Petersen graph}
One of the most famous graphs in graph theory, Petersen graph (see Fig.~\ref{pic4}), 
is SRG with parameter set $(10,3,0,1)$. 
It is the complementary graph to the triangular graph $T(5)$
and it is known to be unique with this parameter set. All of its
equitable partitions have been enumerated in \cite{zv}.
Equation \thetag{eq1} has two solutions: $(a,b)=(2,5)$ and $(a,b)=(5,2)$,
but the second is ruled out by \thetag{eq2}. 
For $(a,b)=(2,5)$, a good equitable partition is determined by a partition
into two disjoint 5-cycles.
The corresponding DSRGs from the $\pi$-join have parameter set 
$(20j+10, 10j+3, 5j+3, 5j, 5j+1)$.
For $j=1$ the parameter set $(30,13,8,5,6)$ is not new, since in \cite{ma} the authors
constructed a DSRG with this parameter set. Their digraph is not vertex-transitive, 
while ours is, thus these graphs are non-isomorphic. 
The first DSRG with parameter set $(50,23,13,10,11)$ has been just recently constructed
in \cite{gy}. For $j\geq 3$ the resulting parameter set is new, according to our best knowledge.

\begin{minipage}{\linewidth}
      \centering
      \begin{minipage}{0.45\linewidth}
          \begin{figure}[H]
              \includegraphics[height=5cm]{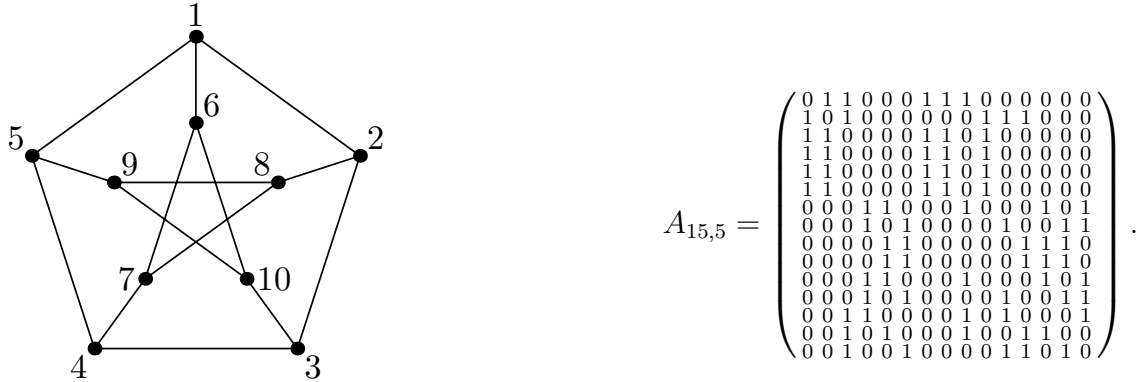}
              \caption{SRG(10,3,0,1).}
              \label{pic4}
          \end{figure}
      \end{minipage}
      \hspace{0.05\linewidth}
      \begin{minipage}{0.45\linewidth}
$$A_{15,5}=\left(\begin{smallmatrix}
0 & 1 & 1 & 0 & 0 & 0 & 1 & 1 & 1 & 0 & 0 & 0 & 0 & 0 & 0 \\
1 & 0 & 1 & 0 & 0 & 0 & 0 & 0 & 0 & 1 & 1 & 1 & 0 & 0 & 0 \\
1 & 1 & 0 & 0 & 0 & 0 & 1 & 1 & 0 & 1 & 0 & 0 & 0 & 0 & 0 \\
1 & 1 & 0 & 0 & 0 & 0 & 1 & 1 & 0 & 1 & 0 & 0 & 0 & 0 & 0 \\
1 & 1 & 0 & 0 & 0 & 0 & 1 & 1 & 0 & 1 & 0 & 0 & 0 & 0 & 0 \\
1 & 1 & 0 & 0 & 0 & 0 & 1 & 1 & 0 & 1 & 0 & 0 & 0 & 0 & 0 \\
0 & 0 & 0 & 1 & 1 & 0 & 0 & 0 & 1 & 0 & 0 & 0 & 1 & 0 & 1 \\
0 & 0 & 0 & 1 & 0 & 1 & 0 & 0 & 0 & 0 & 1 & 0 & 0 & 1 & 1 \\
0 & 0 & 0 & 0 & 1 & 1 & 0 & 0 & 0 & 0 & 0 & 1 & 1 & 1 & 0 \\
0 & 0 & 0 & 0 & 1 & 1 & 0 & 0 & 0 & 0 & 0 & 1 & 1 & 1 & 0 \\
0 & 0 & 0 & 1 & 1 & 0 & 0 & 0 & 1 & 0 & 0 & 0 & 1 & 0 & 1 \\
0 & 0 & 0 & 1 & 0 & 1 & 0 & 0 & 0 & 0 & 1 & 0 & 0 & 1 & 1 \\
0 & 0 & 1 & 1 & 0 & 0 & 0 & 0 & 1 & 0 & 1 & 0 & 0 & 0 & 1 \\
0 & 0 & 1 & 0 & 1 & 0 & 0 & 0 & 1 & 0 & 0 & 1 & 1 & 0 & 0 \\
0 & 0 & 1 & 0 & 0 & 1 & 0 & 0 & 0 & 0 & 1 & 1 & 0 & 1 & 0 \\
\end{smallmatrix}\right).
$$
      \end{minipage}
\end{minipage}

\subsection{DSRG(15,5,2,1,2)}

There are 1292 non-isomorphic DSRGs sharing this parameter set, 
all with very small group of automorphisms (orders less than 6).
We investigated the DSRG corresponding to the adjacency matrix $A_{15,5}$.
Equation \thetag{eq1} has unique solution: $(a,b)=(3,5)$.
A good partition is
$$\{\{1, 2, 7, 8, 10\}, \{3, 4, 9, 11, 15\}, \{5, 6, 12, 13, 14\}\}.$$
The parameter set of the $\pi$-join is 
$(45j+15,15j+5,5j+2,5j+1,5j+2)$, which is new.

\subsection{SRG(16,5,0,2) -- Clebsch graph}

There is a unique SRG with parameters $(16,5,0,2)$, it is the Clebsch graph
(see Fig. \ref{cl}).
Equation \thetag{eq1} has two solutions: $(a,b)=(2,8)$ and $(a,b)=(4,4)$.
The latter solution is ruled out by \thetag{eq2}.
A good partition is $\{\{1, 2, \ldots, 8\}, \{9, 10, \ldots, 16\}\}.$
The parameter set of the $\pi$-join is
$(32j+16,16j+5,8j+5,8j,8j+2)$ which also comes from \cite{gh}.

\begin{minipage}{\linewidth}
      \centering
      \begin{minipage}{0.4\linewidth}
          \begin{figure}[H]
              \includegraphics[height=6cm]{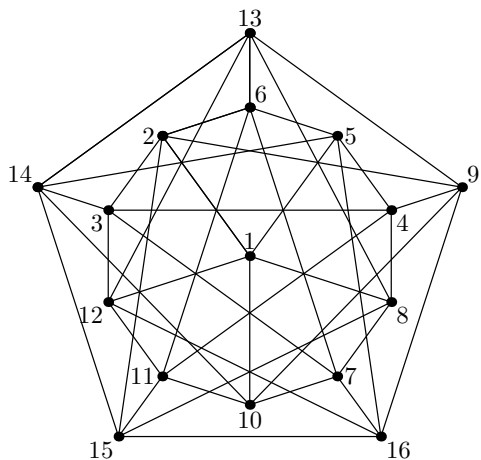}
              \caption{Clebsch graph.}
              \label{cl}  
          \end{figure}
      \end{minipage}
      \hspace{0.05\linewidth}
      \begin{minipage}{0.4\linewidth}
          \begin{figure}[H]
              \includegraphics[height=6cm]{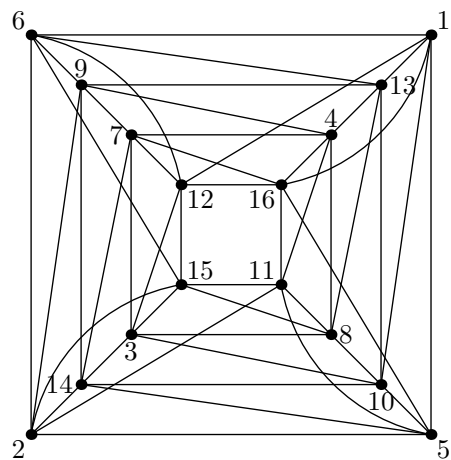}
              \caption{Shrikhande graph.}
              \label{pic7}  
          \end{figure}
      \end{minipage}
\end{minipage}

\subsection{SRG(16,6,2,2)}

There are precisely two non-isomorphic graphs sharing this parameter set:
the Shkrikhande graph (see Fig.~\ref{pic7}) and the lattice square graph
$L_2(4)$ (considered in Section \ref{ss7.6}). 
Using the labeling from the figure, the partitions
$$\{\{1, 2, 5, 6, 9, 10, 13, 14\}, \{3, 4, 7, 8, 11, 12, 15, 16\}\},$$
and 
$$\{\{1, 2, 3, 4\}, \{5, 6, 7, 8\}, \{9, 10, 11, 12\}, \{13, 14, 15, 16\}\}$$
are good partitions of the Shrikhande graph for $(a,b)=(2,8)$ and $(a,b)=(4,4)$, respectively. 
First case lead to the parameter set $(32j+16, 16j+6, 8j+6, 8j+2, 8j+2)$, 
which is known from \cite{gh}, and the second case to the parameter set
$(64j+16,16j+6,4j+6,4j+2,4j+2)$, which is new. 

\subsection{DSRG(16,7,5,4,2)}
There is a unique DSRG on this parameter set, its adjacency matrix 
is $A_{16}$. 
Equation \thetag{eq1} has two solutions: $(a,b)=(2,8)$ and $(a,b)=(4,4)$.
In the case when $(a,b)=(2,8)$ a good partition is:
$$\{\{1, 2, 5, 6, 9, 10, 13, 14\}, \{3, 4, 7, 8, 11, 12, 15, 16\}\}.$$
The resulting parameter sets is $(32j+16,16j+7,8j+5,8j+4,8j+2)$, 
which is known from \cite{du}.
In the case when $(a,b)=(4,4)$ a good partition
$$\{\{1,2,3,4\},\{5,6,7,8\},\{9,10,11,12\},\{13,14,15,16\}\}.$$
The resulting parameter set 
$(64j+16,16j+7,4j+5,4j+4,4j+2)$ is new.

$$A_{16}=\left(\begin{smallmatrix}
0 & 1 & 0 & 0 & 1 & 1 & 0 & 0 & 1 & 1 & 0 & 0 & 1 & 1 & 0 & 0 \\
1 & 0 & 0 & 0 & 1 & 1 & 0 & 0 & 1 & 1 & 0 & 0 & 1 & 1 & 0 & 0 \\
0 & 0 & 0 & 1 & 0 & 0 & 1 & 1 & 0 & 0 & 1 & 1 & 0 & 0 & 1 & 1 \\
0 & 0 & 1 & 0 & 0 & 0 & 1 & 1 & 0 & 0 & 1 & 1 & 0 & 0 & 1 & 1 \\
1 & 1 & 0 & 0 & 0 & 1 & 0 & 0 & 1 & 1 & 0 & 0 & 1 & 1 & 0 & 0 \\
1 & 1 & 0 & 0 & 1 & 0 & 0 & 0 & 1 & 1 & 0 & 0 & 1 & 1 & 0 & 0 \\
0 & 0 & 1 & 1 & 0 & 0 & 0 & 1 & 0 & 0 & 1 & 1 & 0 & 0 & 1 & 1 \\
0 & 0 & 1 & 1 & 0 & 0 & 1 & 0 & 0 & 0 & 1 & 1 & 0 & 0 & 1 & 1 \\
1 & 1 & 0 & 0 & 0 & 0 & 1 & 1 & 0 & 1 & 0 & 0 & 0 & 0 & 1 & 1 \\
1 & 1 & 0 & 0 & 0 & 0 & 1 & 1 & 1 & 0 & 0 & 0 & 0 & 0 & 1 & 1 \\
0 & 0 & 1 & 1 & 1 & 1 & 0 & 0 & 0 & 0 & 0 & 1 & 1 & 1 & 0 & 0 \\
0 & 0 & 1 & 1 & 1 & 1 & 0 & 0 & 0 & 0 & 1 & 0 & 1 & 1 & 0 & 0 \\
0 & 0 & 1 & 1 & 1 & 1 & 0 & 0 & 0 & 0 & 1 & 1 & 0 & 1 & 0 & 0 \\
0 & 0 & 1 & 1 & 1 & 1 & 0 & 0 & 0 & 0 & 1 & 1 & 1 & 0 & 0 & 0 \\
1 & 1 & 0 & 0 & 0 & 0 & 1 & 1 & 1 & 1 & 0 & 0 & 0 & 0 & 0 & 1 \\
1 & 1 & 0 & 0 & 0 & 0 & 1 & 1 & 1 & 1 & 0 & 0 & 0 & 0 & 1 & 0 \\
\end{smallmatrix}\right),\quad
A_{18,7}=\left(\begin{smallmatrix}
0 & 1 & 1 & 0 & 1 & 0 & 0 & 0 & 1 & 1 & 0 & 0 & 0 & 1 & 0 & 0 & 0 & 1 \\
1 & 0 & 1 & 0 & 0 & 1 & 1 & 0 & 0 & 0 & 1 & 0 & 0 & 0 & 1 & 1 & 0 & 0 \\
1 & 1 & 0 & 1 & 0 & 0 & 0 & 1 & 0 & 0 & 0 & 1 & 1 & 0 & 0 & 0 & 1 & 0 \\
0 & 0 & 1 & 0 & 1 & 1 & 0 & 1 & 0 & 0 & 0 & 1 & 1 & 0 & 0 & 0 & 1 & 0 \\
1 & 0 & 0 & 1 & 0 & 1 & 0 & 0 & 1 & 1 & 0 & 0 & 0 & 1 & 0 & 0 & 0 & 1 \\
0 & 1 & 0 & 1 & 1 & 0 & 1 & 0 & 0 & 0 & 1 & 0 & 0 & 0 & 1 & 1 & 0 & 0 \\
0 & 1 & 0 & 0 & 0 & 1 & 0 & 1 & 1 & 0 & 1 & 0 & 0 & 0 & 1 & 1 & 0 & 0 \\
0 & 0 & 1 & 1 & 0 & 0 & 1 & 0 & 1 & 0 & 0 & 1 & 1 & 0 & 0 & 0 & 1 & 0 \\
1 & 0 & 0 & 0 & 1 & 0 & 1 & 1 & 0 & 1 & 0 & 0 & 0 & 1 & 0 & 0 & 0 & 1 \\
1 & 0 & 0 & 0 & 0 & 1 & 0 & 1 & 0 & 0 & 1 & 1 & 0 & 0 & 1 & 0 & 1 & 0 \\
0 & 1 & 0 & 1 & 0 & 0 & 0 & 0 & 1 & 1 & 0 & 1 & 1 & 0 & 0 & 0 & 0 & 1 \\
0 & 0 & 1 & 0 & 1 & 0 & 1 & 0 & 0 & 1 & 1 & 0 & 0 & 1 & 0 & 1 & 0 & 0 \\
0 & 1 & 0 & 1 & 0 & 0 & 0 & 0 & 1 & 0 & 1 & 0 & 0 & 1 & 1 & 0 & 0 & 1 \\
0 & 0 & 1 & 0 & 1 & 0 & 1 & 0 & 0 & 0 & 0 & 1 & 1 & 0 & 1 & 1 & 0 & 0 \\
1 & 0 & 0 & 0 & 0 & 1 & 0 & 1 & 0 & 1 & 0 & 0 & 1 & 1 & 0 & 0 & 1 & 0 \\
0 & 0 & 1 & 0 & 1 & 0 & 1 & 0 & 0 & 0 & 0 & 1 & 0 & 1 & 0 & 0 & 1 & 1 \\
1 & 0 & 0 & 0 & 0 & 1 & 0 & 1 & 0 & 1 & 0 & 0 & 0 & 0 & 1 & 1 & 0 & 1 \\
0 & 1 & 0 & 1 & 0 & 0 & 0 & 0 & 1 & 0 & 1 & 0 & 1 & 0 & 0 & 1 & 1 & 0 \\
\end{smallmatrix}\right).$$

\subsection{DSRG(18,4,3,0,1)}

Duval in his seminal paper presented also a DSRG with parameter set
$(18,4,3,0,1)$. This graph has been known already from other problems in
graph theory, and often is referred to as \emph{Bos\'ak graph} \cite{bk}.
It is known to be unique on this parameter set. 

Though \thetag{eq1} has two solutions, one of them is ruled out by
condition \thetag{eq2}. For $(a,b)=(3,6)$ a good partition is 
given by the partition into three vertex-disjoint induced copies 
of the subgraph isomorphic to the unique DSRG(6,2,1,0,1). 
Hence the existence of a DSRG with parameters 
$(54j+18, 18j+4, 6j+3, 6j, 6j+1)$ is proven. For $j=1$ 
a DSRG with parameter set $(72,22,9,6,7)$ has been recently found 
in \cite{gy}.

\subsection{DSRG(18,7,5,2,3)}

For a DSRG(18,7,5,2,3) with adjacency matrix $A_{18,7}$ we have to consider
two possibilities: $(a,b)=(2,9)$ and $(a,b)=(3,6)$. 
A good partition is $\{\{1,2,\ldots,9\},\{10,\ldots,18\}\}$ in the first case,
and $\{\{1,4,7,10,13,16\},\{2,5,8,11,14,17\},\{3,6,9,12,15,18\}\}$ in the second case.
The resulting parameter sets $(36j+18,18j+7,9j+5,9j+2,9j+3)$ and 
$(54j+18,18j+7,6j+5,6j+2,6j+3)$ are both new, in general. 

\subsection{J\o rgensen's sporadic examples}
\label{ss7.8}

J\o rgensen in \cite{jo} besides a remarkable amount of non-existence results
also constructed a few DSRGs as Cayley graphs. 

\subsubsection{DSRG(21,6,2,1,2)}

Let us consider his DSRG(21,6,2,1,2), which was constructed as a Cayley graph
over the group $G=\langle x, y\,|\, x^3=y^7=1, xy^2=yx\rangle$ with connection set
$S=\{y,y^3,x,xy^2,x^2,x^2y^5\}$. After solving the equation \thetag{eq1} and
counting the quotient matrix we figure out that
$(a,b)=(3,7)$ is a solution and $\lambda+b-k=\mu=2$. 
Therefore the partition $\pi=\{C_0,C_1,C_2\}$ where $C_i=\{x^iy^j\,|\, j=0,1,2,\ldots,6\}$
is a good partition leading to the existence of a family of DSRGs with parameter set
$(63j+21,21j+6,7j+2,7j+1,7j+2)$.

\subsubsection{DSRGs of order 24}
J\o rgensen (\cite{jo}, Theorem 14) constructed three DSRGs on 24 vertices 
as Cayley graphs of $S_4$. Their connection sets are
\[
X_1 = \{ (3, 4), (2, 3), (2, 3, 4), (1, 2)(3, 4), (1, 2, 3, 4), (1, 3, 2), (1, 3, 4, 2), (1, 3, 4)\}, 
\]
\[
X_2 = \{ (3, 4), (2, 3), (1, 2, 3), (1, 2, 4, 3), (1, 3, 2), (1, 3, 4), (1, 3)(2, 4), (1, 4, 3), (1, 4, 2, 3)\},
\]
and
\[
X_3 = \{ (3, 4), (2, 3), (2, 3, 4), (2, 4, 3), (1, 2), (1, 2, 3), (1, 2, 3, 4), (1, 4, 3, 2), (1, 4, 3), (1, 4)\}.
\]
The corresponding Cayley graphs are DSRGs with parameters $(24,8,3,2,3)$, $(24,9,7,2,4)$,
and $(24,10,8,4,4)$, respectively. 

For DSRG(24,8,3,2,3) the equation \thetag{eq1} has only one solution $(a,b)=(3,8)$, 
while for (24,9,7,2,4) and (24,10,8,4,4) there are two solutions: $(2,12)$ and $(3,8)$.

In order to get a good partition for the DSRG(24,8,3,2,3) we run an exhaustive computer search.
It returned us hundreds of such partitions. One of them has cells: \newline
$C_1=\{id, (3,4), (1,2,3), (2,3,4), (1,2,3,4), (1,2,4), (2,4,3), (2,3)\}$, \newline
$C_2=\{(1,3)(2,4), (2,4), (1,3,2,4), (1,3,2), (1,3,4,2), (1,2), (1,4,3,2), (1,2)(3,4)\}$,
\newline 
$C_3=\{(1,4,2), (1,4,3), (1,2,4,3), (1,3), (1,4), (1,3,4), (1,4,2,3), (1,4)(2,3)\}$.

Thus the existence of DSRGs with parameter set
$(72j+24, 24j+8, 8j+3, 8j+2, 8j+3)$ is proved.

For the remaining two DSRGs we have found good partitions combinatorially.

Using conditions on the quotient matrix we immediately get that the partition of $S_4$ into 
even an odd permutations gives us a good partition in the case of DSRG(24,9,7,2,4). 
Thus the existence of a DSRG with parameters $(48j+24,24j+9,12j+7,12j+2,12j+4)$ is proved.
For $(a,b)=(3,8)$ the situation is a little bit more complicated. After some
computer algebra experimentation we concluded, that considering the index 3 subgroup
$H=\langle(1,2)(3,4), (1,2,3,4)\rangle\leq S_4$ the partition $\pi=\{H, H\cdot(2,3),
H\cdot(3,4)\}$ is good. Hence, the existence
of DSRGs with parameter set $(72j+24, 24j+9, 8j+7, 8j+2, 8j+4)$ is proved. 

Let us put permutation $p\in S_4$ into $C_1$ if $1^p\in\{1,2\}$, 
otherwise put $p$ into $C_2$. Partition $\pi=\{C_1,C_2\}$ is a good
partition for DSRG(24,10,8,4,4) with $(a,b)=(2,12)$. This proves the existence of
a DSRG with parameter set $(48j+24, 24j+10, 12j+8, 12j+4, 12j+4)$. 
Considering the index 3 subgroup $H=\langle (1,2,3,4), (1,3)\rangle\leq S_4$ 
the partition $\pi=\{H, H\cdot(1,2), H\cdot(1,4)\}$ is good for DSRG(24,10,8,4,4)
for $(a,b)=(3,8)$, therefore
the existence of DSRG with parameters $(72j+24,24j+10,8j+8,8j+4,8j+4)$ is proved. 

In all the cases, the resulting parameter sets for $j=1$ has been displayed in \cite{ab}
as open.

\subsubsection{DSRG(26,11,7,4,5)}

J\o rgensen (\cite{jo}, Theorem 15) constructed a DSRG(26,11,7,4,5).
He prescribed adjacencies between elements of two sets of size 13. 
For this case \thetag{eq1} says that $(a,b)=(2,13)$ and from the quotient matrix 
it is easy to check that the two original sets define a good partition.  
This confirms the existence of DSRGs with parameters 
$(52j+26,26j+11,13j+7,13j+4,13j+5)$. For $j=1$ the parameter set $(78,37,20,17,18)$
is new according to \cite{ab}.

\subsection{SRG(26,10,3,4)}

According to Paulus \cite{pa} there are precisely 10 SRGs with parameters
(26,10,3,4). The second most symmetric of them has group of automorphisms
of order 39, which has two orbits of size 13 on the vertex set. 
This orbit-partition is good for the unique solution $(a,b)=(2,13)$ of the equation 
\thetag{eq1}, hence the existence of DSRGs with parameters 
$$(52j+26, 26j+10, 13j+10, 13j+3, 13j+4)$$
is proved for all positive integer $j$.
For $j=1$ the resulting parameter set (78,36,23,16,17) is new according to \cite{ab}.

\subsection{SRG(27,10,1,5) -- Complement to the Schl\"afli graph}

The equation \thetag{eq1} gives us a unique solution $(a,b)=(3,9)$. 
For diagonal elements of the quotient matrix of a good partition 
we get $\lambda+b-k=0$, what means that in the basic graph we have
to find a coclique of size 9, or equivalently, a clique of size 9
in the Schl\"afli graph. This graph has been well-studied, and we know
that the maximal clique is of size~6. 
Therefore there is no good partition for our construction. 

\subsection{Sporadic examples by Martinez and Araluze}
\label{ss7.9}

Martinez and Araluze \cite{ma} constructed a handful of DSRGs with previously unknown 
parameter sets with order at most 45 using partial sum families over cyclic groups. 
From their construction immediately follows that the equivalence classes
modulo $a$ (where $a$ is the number of blocks) form a column equitable partition,
and these partitions are not excluded neither by condition \thetag{eq1} nor by 
\thetag{eq2}. 
Moreover, from the blocks it is easy to count the quotient matrix, since it is enough
to count the amount of numbers modulo $a$ in each block for each possibility.

In this way we checked all their 14 graphs constructed from cyclic groups, 
and figured out that six of them have quotient matrix satisfying part (b) in Theorem 1. 
These correspond to parameter sets: 
(27,10,6,3,4), (30,13,8,5,6), (34,14,12,5,6), (34,15,9,6,7),
(39,10,6,1,3) and (40,17,11,6,8).
Consequently, using our $\pi$-join construction we can create 
for any positive integer $j$ a DSRG with parameter set: 
$(81j+27,27j+10,9j+6,9j+3,9j+4)$, $(60j+30,30j+13,15j+8,15j+5,15j+6)$,
$(68j+34,34j+14,17j+12,17j+5,17j+6)$, $(68j+34,34j+15,17j+9,17j+6,17j+7)$, 
$(117j+39,39j+10,13j+6,13j+1,13j+3)$ and $(80j+40,40j+17,20j+11,20j+6,20j+8)$. 

According to \cite{ab} we constructed the first DSRGs with parameters
$(102,48,29,22,23)$, $(102,49,26,23,24)$, $(108,37,15,12,13)$, while the parameter
set $(90,43,23,20,21)$ appears only in this paper, but here already for the third time. 

\subsection{DSRG(36,13,11,2,6)}

The first DSRG with parameter set $(36,13,11,2,6)$ has been recently constructed in \cite{gy}
as a union of classes in a non-commutative association scheme.
The equation \thetag{eq1} has two solutions $(a,b)=(2,18)$ and $(3,12)$. 
The full group of automorphisms of the graph is transitive of order 144,
therefore we can consider its system of imprimitivity.
One can obtain a good partition for the first solution as a union of a 
block of imprimitivity of size 12, and a suitable block of size 6.  
Hence the existence of DSRGs with parameter sets 
$$(72j+36, 36j+13, 18j+11, 18j+2, 18j+6)$$
is confirmed for all positive integer $j$. 
For $j=1$ the parameter set $(108,49,29,20,24)$ is new according to \cite{ab}.

For the case $(a,b)=(3,12)$ a good partition is given by the system
of imprimivity with 3 cells of size 12, where each cell is inducing 
a regular graph of valency 1.
This confirms the existence of DSRGs with parameters 
$$(108j+36,36j+13, 12j+11, 12j+2, 12j+6)$$
for all positive integer $j$.

\subsection{SRG(36,14,4,6)}

There are precisely 180 SRGs (for more information see \cite{bs, sp})
with parameters (36,14,4,6). 
For this parameter set \thetag{eq1} has two solutions:
$(a,b)=(2,18)$ and $(3,12)$. 

Let us consider the second most symmetric graph of them, i.e.
which has group of automorphisms of order 432. 
The imprimitivity system of the group of automorphisms
has blocks of size 18. The corresponding partition fits the
conditions for being good for our construction, therefore 
we can construct for every positive integer $j$ a DSRG 
with parameter set 
$$(72j+36, 36j+14, 18j+14, 18j+4, 18j+6).$$
For $j=1$ the resulting parameter set $(108,50,32,22,24)$ is 
new according to our best knowledge.

If we consider the only imprimitivity system of the group of 
automorphisms with blocks of size 6, then after some experimentation
we can obtain 3 pairs of blocks whose unions results in a 3-cell partition
with cells of size 12 and quotient matrix $2J+6(J-I)$. This proves 
the existence of a DSRG with parameter set 
$$(108j+36, 36j+14, 12j+14, 12j+4, 12j+6)$$
for arbitrary positive integer $j$.  

\subsection{SRG(36,15,6,6)}

It is known that the exact number of SRG sharing the parameter set (36,15,6,6)
is 32548 (for more information see \cite{bs,sp}).
We considered a highly symmetric one of them with automorphism group of order 648.
(In the catalogue of E. Spence it is the last graph on this parameter set.) 
It has a unique system of imprimitivity with blocks of size 9. 
Union of any two of these blocks and the remaining two blocks gives us
a suitable partition for $(a,b)=(2,18)$. 
Hence, the existence of a DSRG with parameter set 
$$(72j+36, 36j+15, 18j+15, 18j+6, 18j+6)$$ is proved for every integer $j$.
For $j=1$ the parameter set $(108,51,33,24,24)$ is new according to \cite{ab}.  

In a similar fashion, using unions of 4 blocks of imprimitivity of size 3 we 
have found a good partition into 3 cells of size 12, therefore the existence
of a DSRG with parameter set 
$$(108j+36, 36j+15, 12j+15, 12j+6, 12j+6)$$  
is proved for all integer $j$.

\subsection{SRG(50,7,0,1) -- Hoffman-Singleton graph}

Just for curiosity we considered one of the most famous SRGs, the Hoffman-Singleton
graph. Solutions of \thetag{eq1} are $(5,10)$ and $(10,5)$, 
however the latter is excluded by \thetag{eq2}. 
The first solution challenges us to consider the partition 
consisting from 5 vertex-disjoint copies of the Petersen graph.
By counting the possible quotient matrix it is not hard to see 
that this is really a good partition, therefore there exists a DSRG 
$(250j+50, 50j+7, 10j+7, 10j, 10j+1)$ for each positive integer $j$.

\section{Summary}
Using our $\pi$-join construction introduced in Section \ref{spi} 
we constructed dozens of infinite families of DSRGs from smaller ones. 
Among the constructed graphs with order less than 110 there are 30 
with new parameter sets according to \cite{ab}.
The new parameter sets are displayed in the Appendix.

\section*{Acknowledgement}
The author gratefully acknowledges the contribution of the Scientific Grant
Agency of the Slovak Republic under the grant 1/0151/15, as well as the constribution
of the Slovak Research and Development Agency under the project APVV 0136-12.
This research was also supported by the Project: Mobility - enhancing research, 
science and education at the Matej Bel University, ITMS code: 26110230082, under the
Operational Program Education cofinanced by the European Social Fund.

\section*{Appendix}
\begin{table}[h]
\begin{center}
\begin{tabular}[h]{|c|c|ccc|r|}
\hline
Parameter set & Parameter set &  &&\\
of the $\pi$-join &  of the basic graph & $a$ & $b$ & $j$ \\
\hline
(40,11,4,3,3) & (8,3,2,1,1) & 4 & 2 & 1  \\
(50,23,13,10,11) & (10,3,0,1) & 2 & 5 & 2  \\
(60,15,4,3,4) & (12,3,1,0,1) & 4 & 3 & 1  \\
(60,20,7,6,7) & (6,2,1,0,1) & 3 & 2 & 3  \\
(60,20,7,6,7) & (15,5,2,1,2) & 3 & 5 & 1  \\  
(70,33,18,15,16) & (10,3,0,1) & 2 & 5 & 3  \\
(72,19,6,5,5) & (8,3,2,1,1) & 4 & 2 & 2  \\
(72,25,11,8,9) & (18,7,5,2,3) & 3 & 6 & 1  \\
(72,33,19,14,16) & (24,9,7,2,4) & 2 & 12 & 1  \\
(72,34,20,16,16) & (24,10,8,4,4) & 2 & 12 & 1  \\
(78,26,9,8,9) & (6,2,1,0,1) & 3 & 2 & 4  \\
(78,36,23,16,17) & (26,10,3,4) & 2 & 13 & 1  \\
(78,37,20,17,18) & (26,11,7,4,5) & 2 & 13 & 1  \\
(80,22,10,6,6) & (16,6,2,2) & 4 & 4 & 1  \\
(80,23,9,8,6) & (16,7,5,4,2) & 4 & 4 & 1  \\
(90,19,5,4,4) & (15,4,2,1,1) & 5 & 3 & 1  \\
(90,31,13,10,11) & (9,4,1,2) & 3 & 3 & 3  \\
(90,43,23,20,21) & (10,3,0,1) & 2 & 5 & 4  \\
(90,43,23,20,21) & (18,7,5,2,3) & 2 & 9 & 2  \\
(90,43,23,20,21) & (30,13,8,5,6) & 2 & 15 & 1  \\
(96,32,11,10,11) & (6,2,1,0,1) & 3 & 2 & 5  \\
(96,32,11,10,11) & (24,8,3,2,3) & 3 & 8 & 1  \\
(96,33,15,10,12) & (24,9,7,2,4) & 3 & 8 & 1  \\
(96,34,16,12,12) & (24,10,8,4,4) & 3 & 8 & 1  \\
(102,48,29,22,23) & (34,14,12,5,6) & 2 & 17 & 1  \\
(102,49,26,23,24) & (34,15,9,6,7) & 2 & 17 & 1  \\
(104,27,8,7,7) & (8,3,2,1,1) & 4 & 2 & 3  \\
(105,35,12,11,12) & (15,5,2,1,2) & 3 & 5 & 2  \\  
(108,27,7,6,7) & (12,3,1,0,1) & 4 & 3 & 2  \\
(108,37,15,12,13) & (27,10,6,3,4) & 3 & 9 & 1  \\
(108,49,29,20,24) & (36,13,11,2,6) & 2 & 18 & 1  \\
(108,50,32,22,24) & (36,14,4,6) & 2 & 18 & 1  \\
(108,51,33,24,24) & (36,15,6,6) & 2 & 18 & 1  \\
(110,53,28,25,26) & (10,3,0,1) & 2 & 5 & 5  \\
\hline
\end{tabular}
\end{center}
\caption{New parameter sets up to $n\leq 110$ appearing in this paper.}
\end{table}

\end{document}